\documentclass[12pt]{article}
\usepackage{amsmath}
\usepackage{amsfonts}

\usepackage{graphicx}
\usepackage{color}
\usepackage{epsfig}\usepackage{psfrag}

\newtheorem{theorem}{Theorem}[section]

\newtheorem{corollary}[theorem]{Corollary}

\newtheorem{definition}[theorem]{Definition}
\newtheorem{example}[theorem]{Example}

\newtheorem{lemma}[theorem]{Lemma}

\newenvironment{proof}[1][Proof]{\textbf{#1.} }{\ \rule{0.5em}{0.5em}}
\renewcommand{\text}{\mbox}

\newcommand{\dom}{\mathop{\rm Dom}}

\newcommand{\graph}{\mathop{\rm Graph}}
\newcommand{\antigraph}{\mathop{\rm Antigraph}}

\newcommand{\R}{\mathbf{R}}

\newcommand{\Z}{\mathbf{Z}}

\begin{document}

\author{Najma Ahmad\thanks{%
Ernst and Young, Toronto Ontario Canada ,  \texttt{%
najma.ahmad@gmail.com}}, Hwa Kil Kim\thanks{%
Courant Institute, New York University, New York NY 10012 USA \texttt{hwakil@cims.nyu.edu}},
and Robert J. McCann\thanks{%
Department of Mathematics, University of Toronto, Toronto Ontario M5S 2E4 Canada, \texttt{%
mccann@math.toronto.edu}}}

\title{Extremal doubly stochastic measures and optimal transportation\thanks{%
It is a pleasure to thank Nassif Ghoussoub and Herbert Kellerer,  who provided early
encouragement in this direction,  and Pierre-Andre Chiappori, Ivar Ekeland, and Lars Nesheim,
whose interest in economic applications fortified our resolve to persist.
We thank Wilfrid Gangbo, Jonathan Korman, and Robert Pego
for fruitful discussions, Nathan Killoran for
useful references, and programs of the Banff International Research Station (2003) and
Mathematical Sciences Research Institute in Berkeley (2005) for stimulating these
developments by bringing us together.
The authors are pleased to acknowledge the support of
Natural Sciences and Engineering Research Council of Canada Grants 217006-03 and -08
and United States National Science Foundation Grant DMS-0354729. \copyright 2009 by the authors.}}
\date{\today}

\begin{abstract}
This article connects the theory of extremal doubly stochastic
measures to the geometry and topology of optimal transportation.

We begin by reviewing an old question (\# 111) of Birkhoff in probability and
statistics \cite{Birkhoff48}, which is to give a necessary and sufficient condition
on the support of a joint probability to guarantee extremality among all
measures which share its marginals. Following work of
Douglas, 
Lindenstrauss, 
and Bene\v{s} and \v{S}t\v{e}p\'an, 
Hestir and Williams~\cite{HestirWilliams95} found a necessary condition which is
nearly sufficient; we relax their subtle
measurability hypotheses separating necessity from sufficiency slightly,
yet demonstrate by example that to be sufficient certainly requires some measurability.
Their condition amounts to the vanishing
of $\gamma$ outside a countable alternating sequence of graphs and antigraphs
in which no two graphs (or two antigraphs) have domains that overlap, and
where the domain of each graph / antigraph in the sequence contains the range
of the succeeding antigraph (respectively, graph).  Such sequences are called
{\em numbered limb systems}. Surprisingly, this characterization can be used
to resolve the uniqueness question for optimal transportation
on manifolds with the topology of the sphere.
\end{abstract}

\maketitle

\section{Introduction}

An $n \times n$ {\em doubly stochastic matrix} refers to a matrix of non-negative entries
whose columns and rows each sum to $1$.  The doubly stochastic
matrices form a convex subset of all $n \times n$ matrices --- in fact a
convex polytope, whose extreme points are in bijective correspondence with the $n!$
permutations on $n$-letters,  according to a theorem of Birkhoff \cite{Birkhoff46} and
von Neumann \cite{vonNeumann53}.
For example, the $3 \times 3$ doubly stochastic matrices,
$$
\left(
\begin{matrix}
s & t & 1-s-t \\
u & v & 1-u-v \\
1-s-u & 1-t-v & s+t+u+v -1 \\
\end{matrix}
\right)
$$
form a 4-dimensional polytope with 6 vertices.
Shortly after proving this characterization,
Birkhoff  \cite[Problem 111]{Birkhoff48} initiated the search for a
infinite-dimensional generalization, thus
stimulating a line of research which remains fruitful even today.

\begin{figure}[h]
\psfragscanon
\centering
\psfrag{o}{$o$}
\epsfig{file=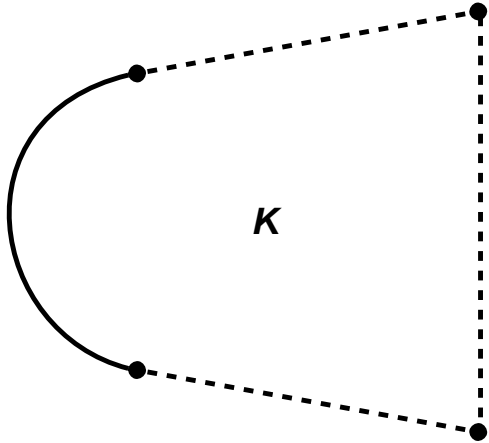, height=6cm}
\caption[Extreme points]
{\label{fig.extreme} Krein-Milman asserts a compact convex set $\mathbf{K}$ can be reconstructed from its extreme
points (denoted here by solid circles $\bullet$ and solid lines $\boldsymbol{-}$)}
\end{figure}

A {\em doubly stochastic measure} on the square refers to a non-negative Borel
probability measure on $[0,1]^2$ whose horizontal and vertical marginals both coincide
with Lebesgue measure $\lambda$ on $[0,1]$.  The set of doubly stochastic measures
forms a convex set we denote by $\Gamma(\lambda,\lambda)$
(which is weak-$*$ compact in the Banach space dual to
continuous functions $C([0,1]^2)$ normed by their suprema $\|\cdot\|_\infty$).
A measure is said to be {\em extremal} in
$\Gamma(\lambda,\lambda)$ if it cannot be decomposed as a convex combination
$\gamma = (1-t) \gamma_0 + t \gamma_1$
with $0<t<1$ and $\gamma_0,\gamma_1 \in \Gamma(\lambda,\lambda)$, except trivially
with $\gamma_0 = \gamma_1$.  Since the Krein-Milman theorem asserts that convex
combinations of extreme points are dense (in any compact convex subset of a topological
vector space, Figure \ref{fig.extreme}),  it is natural to want to characterize the extreme points of
$\Gamma(\lambda,\lambda)$.  Another motivation for such a characterization
is that every continuous linear functional on $\Gamma(\lambda,\lambda)$ is
minimized at an extreme point.  Whether or not this extremum is uniquely
attained can be an interesting question:  in Figure \ref{fig.extreme} the horizontal coordinate is
minimized at a single point but maximized at two extreme points (and along
the segment joining them).

\begin{figure}[h]
\psfragscanon
\centering
\psfrag{o}{$o$}
\epsfig{file=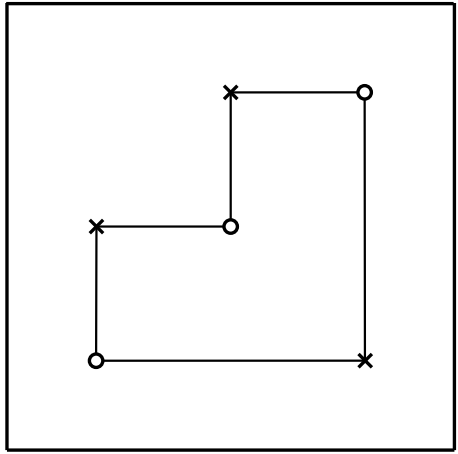, height=6cm}
\caption[Acyclic matrix]
{\label{fig.acyclic} In an {\em acyclic} matrix the product of {\bf x}'s and {\bf o}'s must vanish}
\end{figure}

Motivated by applications like the optimization problem just mentioned,
we prefer to formulate the question in slightly greater generality,  by replacing
the two copies of $([0,1],\lambda)$ with probability spaces
$(X,\mu)$ and $(Y,\nu)$,  where $X$ and $Y$ are each subsets of a complete
separable metric space,  and $\mu$ and $\nu$ are Borel probability measures on
$X$ and $Y$ respectively.  This widens applicability of the answer to this question
without increasing its difficulty.
Letting $\Gamma(\mu,\nu)$ denote the Borel probability measures on $X \times Y$
having $\mu$ and $\nu$ for marginals,  we wish to characterize the extreme points of
the convex set $\Gamma(\mu,\nu)$.  Ideally, as in the finite-dimensional case,
this characterization would be given in terms of some geometrical property of the
support of the measure $\gamma$ in $X \times Y$.  Indeed,  if
$\mu = \sum_{i=1}^m m_i \delta_{x_i}$ and $\nu = \sum_{j=1}^n n_j \delta_{y_j}$
are finite, our problem reduces
to characterizing the extreme points of the convex set ${\cal A}$
of $m \times n$ matrices with prescribed column and row sums:
$$
{\cal A} = \{ a_{ij} \ge 0 \mid m_i = \sum_{j=1}^n a_{ij}, \sum_{i=1}^m a_{ij} = n_j \}.
$$
A matrix $(a_{ij})$ is well-known to be extremal in ${\cal A}$ if and only if it is
{\em acyclic},
meaning for every sequence
$a_{i_1 j_1}, \ldots, a_{i_k j_k}$ of non-zero entries
occupying $k \ge 2$ distinct columns and $k$ distinct rows,
the product $a_{i_1 j_2} \ldots a_{i_{k-1} j_{k}} a_{i_k j_1}$ must vanish
--- see Figure \ref{fig.acyclic} or Denny \cite{Denny80}, where the terminology
{\em aperiodic} is used.  Similarly,  a set $S \subset X \times Y$ is
acyclic if for every $k \ge 2$ distinct points $\{x_1, \ldots, x_k\} \subset X$
and $\{y_1,\ldots,y_k\} \subset Y$,  at least one of the pairs
$(x_1 y_1), (x_1,y_2), (x_2,y_2), \ldots, (x_{k-1},y_k),(x_k, y_k),(x_k,y_1)$
lies outside of $S$.

A functional analytic characterization of extremality was supplied by
Douglas \cite{Douglas64} and by Lindenstrauss \cite{Lindenstrauss65}:
it asserts that $\gamma$ is extremal in $\Gamma(\mu,\nu)$ if and only if
$L^1(X,d\mu) \oplus L^1(Y,d\nu)$ is dense in $L^1(X \times Y,d\gamma)$.
Although this result is a useful starting point,  it is not quite the
characterization we desire for applications,  since it is not easily
expressed in terms of the geometry of the support of $\gamma$.
Significant further progress was made by Bene\v{s}
and \v{S}t\v{e}p\'{a}n, who showed every extremal doubly stochastic
measure vanishes outside some acyclic subset $S \subset X \times Y$ \cite{BenesStepan87}.
Hestir and Williams refined this condition, showing that it
becomes sufficient under an additional Borel measurability hypothesis which,
unfortunately, is not always satisfied \cite{HestirWilliams95}. Some of the subtleties
of the problem were indicated already by Losert's counterexamples \cite%
{Losert82}.  The difficulty of the problem resides partly in
the fact that any geometrical characterization of optimality must be
invariant under arbitrary measure-preserving transformations applied independently
to the horizontal (abscissa) and vertical (ordinate) variables.

In this manuscript we review this line of research, clarifying the nature of
the gap separating necessity from sufficiency and pointing out that
it can be narrowed slightly by replacing the Borel $\sigma$-algebra with suitably
adapted measure-completions. We conclude by describing an application to the question
of uniqueness in optimal transportation,  which is one of the original and most
important examples of an infinite-dimensional program \cite{Kantorovich42},
and appears naturally in applications \cite{RachevRuschendorf98} \cite{Villani09}.
It arises when one wants to use a continuum
of sources to supply a continuum of sinks (modeled by $\mu$ and $\nu$ respectively)
as efficiently as possible.  The question addressed is to identify cost functions $c(x,y)$
on the product space $X \times Y$ whose minimum expected value against measures
in $\Gamma(\mu,\nu)$ is uniquely attained.  When $X$ and $Y$ are differentiable
manifolds and $c \in C^1(X \times Y)$,  to guarantee uniqueness it turns out to
be sufficient that $y_1 \ne y_2$ imply
$x \in X \longrightarrow c(x,y_1)-c(x,y_2)$ has no critical points,
except perhaps for a single global maximum and a single global minimum.
This generalizes to some compact manifolds $X$ a criterion of Gangbo \cite{Gangbo95},
Carlier \cite{Carlier03}, Levin \cite{Levin99} and Ma, Trudinger and Wang \cite{MaTrudingerWang05},
which asserts that the absence of critical points implies uniqueness;
(their condition further implies that almost every source supplies a single sink,
thus solving another transportation problem first posed by Monge \cite{Monge81},
which our condition does not do).
When satisfied,  our criterion implies that
the manifold $X$, if compact, has the topology of the sphere. Uniqueness, however,
remains an interesting open question for compact manifolds which are not topological spheres.
This surprising application was first developed in an economic context
by Chiappori, McCann, and Nesheim \cite{ChiapporiMcCannNesheim08p}.

\section{Measures on graphs are push-forwards}

Before recalling the characterization of interest,  let us develop a bit of notation
in a simpler setting, and a key argument that we shall require.  Impatient or
knowledgeable readers can skim the present section and proceed directly to
the final sections below.

Let $X$ and $Y$ be subsets of complete separable metric spaces,  and fix a
non-negative Borel measure $\mu$ on $X$.
Suppose $f:X \longrightarrow Y$ is {$\mu$-measurable},
meaning $f^{-1}(B)$ is in the $\sigma$-algebra
completion of the Borel subsets of $X$ with respect to the measure $\mu$,
whenever $B$ is relatively Borel in $Y$.
Then a Borel measure on $Y$ is induced,  denoted $f_\#\mu$ and called the
{\em push-forward} of $\mu$ through $f$, and given by
\begin{equation}\label{push-forward}
(f_\#\mu)[B] := \mu[f^{-1}(B)]
\end{equation}
for each Borel $B \subset Y$.  Defining the projections
$\pi^{X_{}}(x,y) = x$ and $\pi^{Y_{}}(x,y)=y$ on $X \times Y$,  this notation permits
the horizontal and vertical marginals of a measure $\gamma \ge 0$ on $X \times Y$
to be expressed as $\pi^{X_{}}_\# \gamma$ and $\pi^{Y_{}}_\# \gamma$ respectively.

The next lemma shows that any measure supported on a graph can be deduced
from its horizontal marginal.  It improves on Lemma 2.4 of \cite{GangboMcCann00}
and various other antecedents,  by using an argument from Villani's
Theorem 5.28 \cite{Villani09}
to extract $\mu$-measurability of $f$ as a conclusion rather that a hypothesis.
As work of, e.g., Hestir and Williams \cite{HestirWilliams95} implies,
although measures on graphs are extremal in $\Gamma (\mu ,\nu )$, the converse is
far from being true; this peculiarity is an inevitable consequence of the infinite
divisibility of $(X,\mu)$.

\begin{lemma}[Measures on graphs are push-forwards]
\label{pure implies unique} Let $X_{}$ and $Y_{}$ be subsets of
complete separable metric spaces, and $\gamma \geq 0$ a $\sigma $-finite
Borel measure on the product space $X_{}\times Y_{}$. Denote the horizontal
marginal of $\gamma $ by $\mu _{}:=\pi _{\#}^{X_{}}\gamma $. If $\gamma $
vanishes outside the graph of $f:X_{}\longrightarrow Y_{}$, meaning
$\{(x,y)\in X_{}\times Y_{}\mid y\neq f(x)\}$ has zero outer measure, then $%
f$ is $\mu _{}$-measurable and $\gamma =(id_{X_{}}\times f)_{\#}\mu _{}$.
\end{lemma}

\begin{proof}
Since outer-measure is subadditive, it costs no generality to assume the
subsets $X_{}$ and $Y_{}$ are in fact complete and separable,
by extending $\gamma$ in the obvious (minimal) way. Any $%
\sigma $-finite Borel measure $\gamma$ is regular and $\sigma$-compact on a
complete separable metric space; e.g.\ p.~255 of \cite{Dudley02} or Theorem
I-55 of \cite{VillaniAnalysis}.
Since $\gamma$ vanishes outside $\graph (f) := \{(x,f(x)) \mid x \in X_{}\}$, there is an increasing sequence
of compact sets $K_i \subset K_{i+1} \subset \graph(f)$ whose
union $K_\infty = \lim_{i \to \infty} K_i$ contains the full mass of $\gamma$%
. Compactness of $K_i \subset \graph (f)$ implies continuity of $%
f$ on the compact projection $X_i := \pi^X(K_i)$. Thus the restriction $%
f_\infty$ of $f$ to $X_\infty := \pi^X(K_\infty)$ is a Borel map whose graph
$K_\infty = \graph(f_\infty)$ is a $\sigma$-compact set of full
measure for $\gamma$. We now verify that $\gamma$ and $(id_{X_\infty} \times
f_\infty)_\# \mu_{}$ assign the same mass to each Borel rectangle $U \times V
\subset X_{} \times Y_{}$. Since $(U \times V) \cap \mathop{\rm Graph}%
(f_\infty) = ((U \cap f_\infty^{-1}(V)) \times Y_{}) \cap \mathop{\rm Graph}%
(f_\infty)$ we find
\begin{eqnarray*}
\gamma(U \times V) &=& \gamma((U \cap f_\infty^{-1}(V)) \times Y_{}) \\
&=& \mu_{}(U \cap f_\infty^{-1}(V)),
\end{eqnarray*}
proving $\gamma = (id_{X_\infty} \times f_\infty)_\# \mu_{}$. Taking $U= X_{}
\setminus X_\infty$ and $V = Y_{}$ shows $X_{} \setminus X_\infty$ is $\mu_{}$%
-negligible. Since $id_{X_{}} \times f$ differs from the Borel map $%
id_{X_\infty} \times f_\infty$ only on the $\mu_{}$-negligible complement of
the $\sigma$-compact set $X_\infty$, we conclude $f$ is $\mu_{}$-measurable
and $\gamma = (id_{X_{}} \times f)_\# \mu_{}$ as desired.
\end{proof}

The preceding lemma shows that any measure concentrated on a graph is
uniquely determined by its marginals; $\gamma$ is therefore extremal in
$\Gamma(\pi^{X_{}}_\#\gamma,\pi^{Y_{}}_\#\gamma)$.  As the results of the
next section show,  the converse is far from being true.

\section{Numbered limb systems and extremality}

In this section we adapt Hestir and Williams \cite{HestirWilliams95}
notion of a {\em numbered limb system} to $X \times Y$.  Using the axiom of choice,
Hestir and Williams deduced from the acyclicity condition of Bene\v{s} and \v{S}t\v{e}p\'an
\cite{BenesStepan87} that each extremal doubly stochastic measure vanishes
outside some numbered limb system.  Conversely,  they showed that vanishing
outside a number limb system is sufficient to guarantee extremality of a
doubly stochastic measure,  provided the graphs (and antigraphs) comprising
the system are Borel subsets of the square.
Our main theorem gives a new proof of this converse in the more general setting
of subsets $X \times Y$ of complete separable metric spaces, and under a slightly
weaker measurability hypothesis on the graphs and antigraphs.
A simple example shows that some measurability hypothesis is nevertheless
required. In the next section, we shall see how this converse is
germane to the question of uniqueness in optimal transportation.

Given a map $f:D\longrightarrow Y$ on $D\subset X$, we denote its graph,
domain, range, and the graph of its (multivalued) inverse by
\begin{eqnarray*}
\mathop{\rm Graph}(f):= &\{(x,f(x))\mid x\in D\},& \\
\mathop{\rm Dom}f:= &\pi ^{X}(\mathop{\rm Graph}(f))&=D, \\
\mathop{\rm Ran}f:= &&\phantom{=}\pi ^{Y}(\mathop{\rm Graph}(f)), \\
\mathop{\rm Antigraph}(f):= &\{(f(x),x)\mid x\in \mathop{\rm Dom}f\}&\subset
Y\times X.
\end{eqnarray*}%
More typically, we will be interested in the $\mathop{\rm Antigraph}%
(g)\subset X\times Y$ of a map $g:D\subset Y\longrightarrow X$.

\begin{figure}[h]
\psfragscanon
\centering
\psfrag{o}{$o$}
\epsfig{file=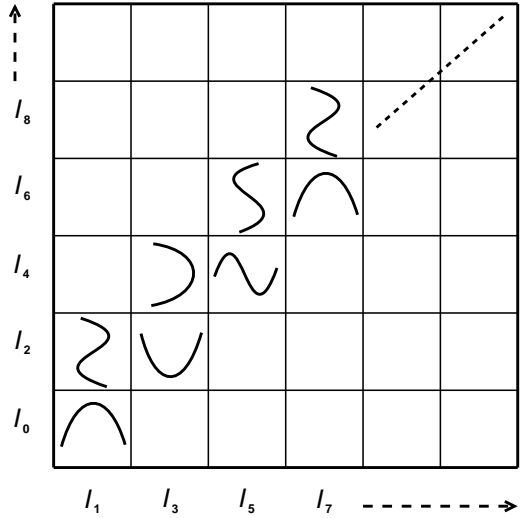, height=7cm}
\caption[Numbered limb system]
{\label{fig.numbered_limb}  \centering The subsets $I_{k}$ need not be connected; in this
numbered limb system they are represented as connected sets for visual convenience only.}
\end{figure}

\begin{definition}[Numbered limb system]
\label{numbered limb system} Let $X_{}$ and $Y_{}$ be Borel subsets of
complete separable metric spaces. A relation $S \subset X_{} \times Y_{}$ is a
\emph{numbered limb system} if there is a countable disjoint decomposition
of $X_{} = \cup_{i=0}^\infty I_{2i+1}$ and of $Y_{} = \cup_{i=0}^\infty I_{2i}$
with a sequence of maps $f_{2i}: \mathop{\rm Dom}(f_{2i}) \subset Y_{}
\longrightarrow X$ and $f_{2i+1}: \mathop{\rm Dom}(f_{2i+1}) \subset X
\longrightarrow Y$ such that $S = \cup_{i=1}^\infty \mathop{\rm Graph}%
(f_{2i-1}) \cup \mathop{\rm Antigraph}(f_{2i})$, with
$\mathop{\rm Dom}(f_{k}) \cup \mathop{\rm Ran}(f_{k+1}) \subset I_{k}$ for each
$k \ge 0$.
The system has (at most) $N$ limbs if $\mathop{\rm Dom}(f_k) = \emptyset$
for all $k>N$. 
\end{definition}

Notice the map $f_0$ is irrelevant to this definition though $I_0$ is not;
we may always take $%
\mathop{\rm Dom}(f_0) = \emptyset$, but require $\mathop{\rm Ran}(f_1)
\subset I_0$. The point is the following theorem and its corollary, which
extends and relaxes the result proved by Hestir and Williams for Lebesgue
measure $\mu_{}=\nu_{} = \lambda$ on the interval $X_{}=Y_{}=[0,1]$.
In it, $\Gamma(\mu,\nu)$ denotes the set of non-negative Borel measures
on $X \times Y$ having $\mu = \pi^X_\# \gamma$ and $\nu = \pi^Y_\# \gamma$
for marginals.
As in the preceding lemma, we say $\gamma$ {\em vanishes}
outside of $S \subset X \times Y$ if $\gamma$ assigns zero outer measure
to the complement of $S$ in $X \times Y$.

\begin{theorem}[Numbered limb systems yield unique correlations]
\label{HestirWilliams}\ \\Let $X_{}$ and $Y_{}$ be subsets of complete
separable metric spaces, equipped with $\sigma$-finite Borel measures
$\mu$ on $X$ and $\nu$ on $Y$.
Suppose there is a numbered limb system $S = \cup_{i=1}^\infty
\mathop{\rm Graph}(f_{2i-1}) \cup \mathop{\rm Antigraph}(f_{2i})$
with the property that $\graph(f_{2i-1})$ and $\antigraph(f_{2i})$
are $\gamma$-measurable subsets of $X \times Y$ for each $i \ge 1$
and for every $\gamma \in \Gamma(\mu,\nu)$ vanishing outside of $S$.
If the system has finitely many limbs or $%
\mu[X]<\infty$, then at most one $\gamma \in \Gamma(\mu,\nu)$
vanishes outside of $S$. If such a measure exists, it is given by
$\gamma = \sum_{k=1}^\infty \gamma_k$ where
\begin{eqnarray}  \label{alternating representation}
\gamma_{2i-1} = (id_{X_{}} \times f_{2i-1})_\# \eta_{2i-1},
&& \gamma_{2i} = (f_{2i} \times id_{Y_{}})_\# \eta_{2i}, \\
 \eta_{2i-1} = \Big(\mu - \pi^{X}_\#\gamma_{2i}\big)   \Big|_{\dom f_{2i-1}} ,
&& \eta_{2i} = \Big(\nu - \pi^{Y}_\#\gamma_{2i+1} \Big) \Big|_{\dom f_{2i}}.
\label{alternating marginals}
\end{eqnarray}
Here $f_k$ is measurable with respect to the $\eta_{k}$ completion of the
Borel $\sigma$-algebra. If the system has $N<\infty$ limbs, $\gamma_k=0$ for
$k > N$, and $\eta_k$ and $\gamma_k$ can be computed recursively from the
formulae above starting from $k=N$.
\end{theorem}

\begin{proof}
Let $S=\cup _{i=1}^{\infty }\mathop{\rm Graph}(f_{2i-1})\cup
\mathop{\rm
Antigraph}(f_{2i})$ be a numbered limb system whose complement has
zero outer measure for some $\sigma$-finite measure $0 \le \gamma \in \Gamma(\mu,\nu)$.
This means that $I_{k} \supset \mathop{\rm Dom} f_{k}$
gives a disjoint decomposition of $X_{}=\cup _{i=0}^{\infty }I_{2i+1}$
and of $Y_{}=\cup _{i=0}^{\infty }I_{2i}$, and that $\mathop{\rm Ran}%
(f_{k})\subset I_{k-1}$ for each $k\geq 1$.
Assume moreover,  that $\graph(f_{2i})$ and $\antigraph(f_{2i-1})$
are $\gamma$-measurable for each $i \ge 1$.  We wish to show $\gamma$ is
uniquely determined by $\mu$, $\nu$ and $S$.

The graphs $%
\mathop{\rm Graph}(f_{2i-1})$ are disjoint since their domains $I_{2i-1}$ are
disjoint, and the antigraphs $\mathop{\rm Antigraph}(f_{2i})$ are disjoint
since their domains $I_{2i}$ are. Moreover, $\mathop{\rm Graph}(f_{2i-1})$
is disjoint from $\mathop{\rm Antigraph}(f_{2j})$ for all $i,j\geq 1$: $%
\mathop{\rm Ran}(f_{2i-1})\subset I_{2i-2}$ prevents
$\mathop{\rm Graph}(f_{2i-1})$ from intersecting
$\mathop{\rm Antigraph}(f_{2j-2})$ unless $%
j=i$ since the domains $I_{2j-2}$ are disjoint, and
$\mathop{\rm Graph}(f_{2i-1})$ cannot intersect $\mathop{\rm Antigraph}(f_{2i-2})$ since $%
\mathop{\rm Dom}(f_{2i-1}) \subset I_{2i-1}$ is disjoint from $\mathop{\rm Ran}%
(f_{2i-2})\subset I_{2i-3}$.

Let $\gamma_k$ denote the restriction of $\gamma$ to $\antigraph(f_k)$
for $k$ even and to $\graph(f_k)$ for $k$ odd.  Then $\gamma = \sum \gamma_{k}$
by our measurability hypothesis,  and $\gamma_k$ restricts to a Borel
measure on $X \times \dom f_k$ if $k$ is even,  and on $\dom f_k \times Y$
if $k$ odd.
Defining the marginal projections $\mu _{k}=\pi _{\#}^{X_{}}\gamma _{k}$ and $%
\nu _{k}=\pi _{\#}^{Y_{}}\gamma _{k}$, setting
$\eta _{k}=\nu _{k}$ if $k$ even and $\eta_{k}=\mu _{k}$ if $k$ odd
yields (\ref{alternating representation}) and the
$\eta _{k}$-measurability of $f_{k}$ immediately from Lemma \ref{pure implies
unique}.
Since $\nu _{2i}$ vanishes outside $\mathop{\rm Dom}f_{2i}$, from
$\nu = \sum_{k=1}^\infty \nu_k$
we derive $\nu _{2i}=(\nu -\sum_{k \neq 2i}\nu _{k})|_{\mathop{\rm Dom}f_{2i}}$. For $k$
even, $\nu _{k}$ vanishes outside $\mathop{\rm Dom}f_{k} \subset I_k$, while for $k$
odd, $\nu _{k}$ vanishes outside $\mathop{\rm Ran}f_{k}\subset I_{k-1}$,
which is disjoint from $\mathop{\rm Dom}f_{2i}$ unless $%
k=2i+1$. Thus $\eta _{2i}=(\nu -\nu _{2i+1})|_{\mathop{\rm Dom}f_{2i}}$. The
formula (\ref{alternating marginals}) for $\eta _{2i-1}$ follows from
similar considerations.

It remains to show the representation (\ref{alternating representation})--(%
\ref{alternating marginals}) specifies $(\gamma _{k},\eta _{k})$ uniquely
for all $k\geq 1$, and hence determines $\gamma =\sum \gamma _{k}$ uniquely.
If the system has $N<\infty $ limbs, $I_{k}=\emptyset $ for $k>N$ and hence $%
\gamma _{k}=0$. We can compute $\eta _{k}$ and $\gamma _{k}$ starting with $%
k=N$, and then recursively from the formulae above for $k=N-1,N-2,\ldots ,1$%
, so the formulae represent $\gamma $ uniquely. If instead $S$ has countably
many limbs, suppose there are two finite Borel measures $\gamma $ and $\bar{%
\gamma}$ vanishing outside of $S$ and having the same marginals $\mu $ and $\nu $.
For each $k \ge 1$,  recall that
\begin{equation*}
K _{k}:=\left\{
\begin{array}{lc}
\graph(f_k) & k\ \mathrm{odd,} \\
\antigraph(f_k) & k\ \mathrm{even,}%
\end{array}%
\right.
\end{equation*}%
is measurable with respect to both $\gamma$ and $\bar \gamma$.
Given $\epsilon >0$, take $N$ large enough so that both $\gamma $ and $\bar{\gamma}$
assign mass less than $\epsilon $ to $\cup _{k=N}^{\infty }K_{k}$.
Set $\gamma_k = \gamma|_{K_k}$ and $\bar{\gamma}_{k}=\bar{\gamma}|_{K_{k}}$
and denote their marginals by $(\mu_k,\nu_k) = (\pi^X_\# \gamma_k,\pi^Y_\# \gamma_k)$
and $(\bar \mu_k,\bar \nu_k) = (\pi^X_\# \bar \gamma_k,\pi^Y_\# \bar \gamma_k)$.
Observe that both
$\gamma ^{\epsilon}:=\sum_{k=1}^{N}\gamma _{k}$ and
$\bar{\gamma}^{\epsilon }:=\sum_{k=1}^{N}\bar{\gamma}_{k}$
are concentrated on the same numbered limb system; it has
finitely many limbs, and the differences
$\delta \mu^{\epsilon }= \sum_{k=1}^N (\bar \mu_k - \mu_k)$ and
$\delta \nu^{\epsilon }= \sum_{k=1}^N (\bar \nu_k - \nu_k)$ between the marginals of
$\gamma^\epsilon$ 
and $\bar \gamma^\epsilon$ 
have total variation at most $2\epsilon $. Since the
$\delta \mu _{2i-1} = \bar \mu_{2i-1} - \mu_{2i-1}$ are mutually singular,
as are the $\delta \nu _{2i} = \bar \nu_{2i} - \nu_{2i}$,
we find the sum of the total variations of
\begin{equation*}
\delta \eta _{k}:=\left\{
\begin{array}{cc}
\bar{\mu}_{k}-\mu _{k} & k\ \mathrm{odd,} \\
\bar{\nu}_{k}-\nu _{k} & k\ \mathrm{even,}%
\end{array}%
\right.
\end{equation*}%
is bounded: $\sum_{k=1}^{N}\Vert \delta \eta _{k}\Vert _{TV(\mathop{\rm Dom}%
f_{k})}<4\epsilon $. Using (\ref{alternating representation}) to derive
\begin{eqnarray*}
\Vert \bar{\gamma}_{k}-\gamma _{k}\Vert
_{TV(X_{}\times Y_{})} &=&\left\{
\begin{array}{cc}
\Vert (id_{X_{}}\times f_{k})_{\#}\delta \eta _{k}\Vert _{TV(X_{}\times
Y_{})} & k\ \mathrm{odd,} \\
\Vert (f_{k}\times id_{Y_{}})_{\#}\delta \eta _{k}\Vert _{TV(X_{}\times
Y_{})} & k\ \mathrm{even,}%
\end{array}%
\right. \\
&=&\Vert \delta \eta _{k}\Vert _{TV(\mathop{\rm Dom}f_{k})}
\end{eqnarray*}%
%
%
%
%
%
%
%
%
%
%
%
%
%
%
%
%
and summing on $k$ yields
$\Vert \bar{\gamma}^{\epsilon }-\gamma ^{\epsilon }\Vert _{TV(X_{}\times
Y_{})}<4\epsilon $. Since $\gamma ^{\epsilon }\rightarrow \gamma $ and $%
\bar{\gamma}^{\epsilon }\rightarrow \bar{\gamma}$ as $\epsilon \rightarrow 0$%
, we conclude $\bar{\gamma}=\gamma $ to complete the uniqueness proof.
\end{proof}

As in Hestir and Williams \cite{HestirWilliams95},  the uniqueness theorem
above implies extremality as an immediate consequence.

\begin{corollary}[Sufficient condition for extremality]
\label{C:support characterization} Let $X$ and $Y$ be subsets of complete separable
metric spaces,  equipped with $\sigma$-finite Borel measures $\mu$ on $X$ and $\nu$ on $Y$.
Suppose there is a numbered limb system $S = \cup_{i=1}^\infty
\mathop{\rm Graph}(f_{2i-1}) \cup \mathop{\rm Antigraph}(f_{2i})$
with the property that $\graph(f_{2i-1})$ and $\antigraph(f_{2i})$
are $\gamma$-measurable subsets of $X \times Y$ for each $i \ge 1$,
for every $\gamma \in \Gamma(\mu,\nu)$ vanishing outside of $S$.
If the system has finitely many limbs or $%
\mu[X]<\infty$, then any measure $\gamma \in \Gamma(\mu,\nu)$
vanishing outside of $S$ is extremal in the convex set $\Gamma(\mu,\nu)$.
\end{corollary}

\begin{proof}
Suppose a measure $\gamma \in \Gamma(\mu,\nu)$ vanishes outside a numbered limb
system $S$ satisfying the hypotheses of the corollary.
If $\gamma =(1-t)\gamma
_{0}+t\gamma _{1}$ with $\gamma _{0},\gamma _{1}\in \Gamma (\mu ,\nu )$ and $%
0<t<1$, then $\gamma \geq \gamma _{0}$ and $\gamma \geq \gamma _{1}$, so both
$\gamma _{0}$ and $\gamma _{1}$ vanish outside of $S$.
According to Theorem \ref{HestirWilliams}, they are uniquely determined by
$S $ and their marginals, hence $\gamma _{0}=\gamma _{1}$ to establish the corollary.
\end{proof}

The following example confirms that a measurability gap still remains between
the necessary and sufficient conditions for extremality.
It is a close variation on the standard example of a non-Lebesgue measurable
set from real analysis.  Together with the lemma and theorem preceding,  this
example makes clear that measurability is required only to allow the graphs
to be separated from each other and from the antigraphs in an additive way.

\begin{example}[An acyclic set supporting non-extremal measures]\ \\
Let $\lambda$ denote Lebesgue measure and define the maps $f_0(x) = x$ and
$f_1(x) = x + \sqrt 2$ (mod 1) on the unit interval $X=Y=[0,1]$.
Notice $\graph(f_i) \subset [0,1]^2$ supports the doubly stochastic measure
$\gamma_i = (id \times f_i)_\#\lambda$ for $i=0$ and $i=1$;  (both measures
are extremal in $\Gamma(\lambda,\lambda)$ by Corollary \ref{C:support characterization}).
Irrationality of $\sqrt 2$ implies $S = \graph(f_0) \cup \graph(f_1)$ is an acyclic set,
hence can be expressed as a numbered limb system according to Hestir and
Williams \cite{HestirWilliams95}. On the other hand,  there are doubly stochastic
measures such as $\gamma:= \frac{1}{2}(\gamma_0+\gamma_1)$ which vanish outside of $S$ but
which are manifestly not extremal.
\end{example}

\section{Uniqueness of optimal transportation}

In this section we illustrate the significance of the foregoing results
by applying them to the uniqueness question for optimal transportation on manifolds.
Given subsets $X$ and $Y$ of complete separable metric spaces equipped with Borel
probability measures,  representing the distributions  $\mu$
of production on $X$ and $\nu$ of consumption on $Y$,  the
Kantorovich-Koopmans \cite{Kantorovich42} \cite{Koopmans49}
transportation problem is to find $\bar \gamma \in\Gamma(\mu,\nu)$ correlating
production with consumption so as to minimize the expected transportation cost
\begin{equation}\label{MKPa}
\inf_{\gamma \in \Gamma(\mu,\nu)} \int_{X \times Y} c(x,y) d\gamma(x,y)
\end{equation}
against some continuous function $c \in C(X \times Y)$.  Hereafter we shall be solely
concerned with the case in which $X$ is a differentiable manifold,  $\mu$ is
absolutely continuous with respect to coordinates on $X$,  and the cost function
$c \in C^1(X \times Y)$ is differentiable with local control on the magnitude of its
$x$-derivative $d_x c(x,y)$ uniformly in $y$;  for convenience we also suppose $Y$ to be a differentiable
manifold and $c$ is bounded,  though though this is not really necessary:
substantially weaker assumptions also suffice \cite{ChiapporiMcCannNesheim08p}.

In this setting one immediately asks whether the infimum (\ref{MKPa}) is
uniquely attained.  Since attainment is evident,  the question here is uniqueness.
If $c$ satisfies a {\em twist} condition, meaning
$x \in X \longrightarrow c(x,y_1) - c(x,y_2)$ has no critical points for
$y_1 \ne y_2 \in Y$,  then not only is the minimizing $\gamma$ unique,  but its
mass concentrates entirely on the graph of a single map $f_1:X \longrightarrow Y$
(a numbered limb system with one limb), thus solving a form of the transportation
problem posed earlier by Monge \cite{Monge81} \cite{Kantorovich48}.
This was proved in comparable generality by Gangbo \cite{Gangbo95}, Carlier \cite{Carlier03}, 
Levin \cite{Levin99}, 
and Ma, Trudinger and Wang \cite{MaTrudingerWang05},  building on the more specific
examples of strictly convex cost functions $c(x,y) = h(x-y)$ in $X=Y=\R^n$
analyzed by Caffarelli \cite{Caffarelli96} and
Gangbo and McCann 
\cite{GangboMcCann96},
and in case $h(x)=|x|^2$ by Abdellaoui and Heinich, Brenier,
Cuesta-Albertos, Matran, and Tuero-Diaz, Cullen and Purser, Knott and Smith,
and R\"uschendorf and Rachev; see \cite{Brenier91} \cite{GangboMcCann96} \cite{Villani09}.
Adding further restrictions beyond this twist hypothesis allowed
Ma, Trudinger, Wang, and later Loeper, to develop
a regularity theory for the map $f_1:X \longrightarrow Y$,
embracing Delanoe, Caffarelli and Urbas' results for the quadratic cost,
Gangbo and McCann's for its restriction to to convex surfaces, and
Wang's the reflector antenna design, which involves the restriction of
$c(x,y)= -\log |x-y|$ to the sphere; references may be found in
\cite{KimMcCann07p} \cite{Villani09}.
Unfortunately,  the twist hypothesis,  also known as a generalized Spence-Mirrlees
condition in the economic literature,
cannot be satisfied for smooth costs $c$ on compact manifolds $X \times Y$,
and apart from the result we are about to discuss there are
no general theorems which
guarantee uniqueness of minimizer to (\ref{MKPa}) in this setting.  With this in mind,
let us state our main theorem,  a version of which was established
in a more complicated economic setting by Chiappori, Nesheim, and
McCann \cite{ChiapporiMcCannNesheim08p}.  We expect the simpler formulation
and argument given below to prove more interesting and accessible to
a mathematical readership.


\begin{theorem}[Uniqueness of optimal transport on manifolds]
\ \\Let $X$ and $Y$ be complete separable manifolds equipped with Borel probability
measures $\mu$ on $X$ and $\nu$ on $Y$.  Let $c \in C^1(X \times Y)$ be a bounded
cost function such that for each $y_1 \ne y_2 \in Y$, the map
\begin{equation}\label{subtwist}
x \in X \longrightarrow c(x,y_1) - c(x,y_2)
\end{equation}
has no critical points,  save at most one global minimum and at most one global maximum.
Assume $d_x c(x,y)$ is locally bounded in $x$,  uniformly in $Y$.
If $\mu$ is absolutely continuous with respect to coordinate measure on $X$,
then the minimum (\ref{MKPa}) is uniquely attained;  moreover,  the minimizer
$\gamma \in \Gamma(\mu,\nu)$ vanishes outside a numbered limb system having at most two
limbs.
\end{theorem}

\begin{proof}
Here we give only the proof that there is a numbered limb system having at most
two limbs,  outside of which the mass of all minimizers $\gamma$ vanishes.
A detailed argument confirming the plausible fact that the graphs of these limbs are
Borel subsets of $X \times Y$ can be found in \cite{ChiapporiMcCannNesheim08p}.
Uniqueness of $\gamma$ then follows from Theorem \ref{HestirWilliams}.

By linear programming duality due to Kantorovich and Koopmans in this context,
it is well-known \cite{Villani09}
that there exist potentials $q \in L^1(X,d\mu)$ and $r \in L^1(Y,d\nu)$ with
\begin{equation}\label{c-transform}
q(x) = \inf_{y \in Y} c(x,y) - r(y)
\end{equation}
such that
\begin{equation}\label{duality}
\inf_{\gamma \in \Gamma(\mu,\nu)} \int_{X \times Y} c(x,y) d\gamma(x,y) = \int_X q(x) d\mu(x) + \int_Y r(y) d\nu(y).
\end{equation}
From (\ref{c-transform}) we see
\begin{equation}\label{zero set}
c(x,y) - q(x) - r(y) \ge 0,
\end{equation}
while (\ref{duality}) implies any minimizer $\gamma \in \Gamma(\mu,\nu)$
vanishes outside the zero set $Z \subset X \times Y$ of the non-negative
function appearing in (\ref{zero set}).  It remains
to show this set $Z$ is contained in a numbered limb system consisting of
at most two limbs (apart from a $\mu \otimes \nu$ negligible set).

From (\ref{c-transform}), $q$ is locally Lipschitz, since
$d_x c(x,y)$ is controlled locally in $x$, independently of $y \in Y$.
Rademacher's theorem therefore combines with absolute continuity of
$\mu$ to imply $q$ is differentiable $\mu$-almost everywhere;  we can
safely ignore any points in $X$ where differentiability of $q$ fails,
since they constitute a set of zero volume: $\gamma[\dom Dq \times Y]
= \mu[\dom Dq]=1$.
Taking $x_0 \in \dom Dq$,  suppose $(x_0,y_1)$ and $(x_0,y_2)$ both lie in $Z$,
hence saturate the inequality (\ref{zero set}).
Then $d_x c(x_0,y_1) = Dq(x_0) = d_x c(x_0,y_2)$.
In case the cost is twisted,  meaning
(\ref{subtwist}) has no critical points,  we conclude $y_1=y_2$ hence
$Z \cap (\dom Dq \times Y)$ is contained in a graph.
This completes the proofs by Gangbo, Carlier, and Ma-Trudinger-Wang,
of existence (and uniqueness) of a solution $y_1=f_1(x_0)$ to Monge's problem,
pairing almost every $x_0 \in X$ with a single $y_1 \in Y$.
Notice uniqueness follows from Lemma \ref{pure implies unique} without further
measurability assumptions.

In the present setting,  however,  we only know that $x_0$ must be
a global minimum or global maximum of the function (\ref{subtwist}).  Exchanging
$y_1$ with $y_2$ if necessary yields
\begin{equation}\label{keep away}
q(x) \le c(x,y_1) - r(y_1)  \le c(x,y_2) - r(y_2)
\end{equation}
for all $x \in X$,  the second inequality being strict unless $x=x_0$,
in which case both inequalities are saturated.  Strictness of inequality
(\ref{keep away}) implies $(x,y_2) \not\in Z$ unless $x=x_0$.  In other words,
$(x,y_2) \in Z$ lies on the antigraph of a function $f_2(y_2) = x_0$
well-defined at $y_2$.  There may or may not be a point
$y_0 \in Y$ different from $y_1$ such that
$$
q(x) \le c(x,y_0) - r(y_0) \le c(x,y_1) - r(y_1)
$$
for all $x \in X$.  If such a point $y_0$ exists,  then
$(x_0,y_1) \in \antigraph(f_2)$ as above.  If no such $y_0$ exists,
setting $f_1(x_0) := y_1$ yields
$Z \cap (\dom Dq \times Y) \subset \graph(f_1) \cup \antigraph(f_2)$.
Since the range of $f_1$ is disjoint from the domain of $f_2$,  this completes
the proof that --- up to $\gamma$-negligible sets ---
$Z$ lies in a numbered limb system with at most two limbs, as desired.
\end{proof}

Let us conclude by recalling an example of an extremal doubly stochastic measure
which does not lie on the graph of a single map,  drawn from
work of Gangbo and McCann \cite{GangboMcCann00} and Ahmad \cite{Ahmad04}
on optimal transportation,
and developed in an economic context by Chiappori, McCann, and
Nesheim~\cite{ChiapporiMcCannNesheim08p}.
Other examples may be found in the work of Seethoff and Shiflett~\cite{SeethoffShiflett78},
Losert~\cite{Losert82}, Hestir and Williams~\cite{HestirWilliams95},
Gangbo and McCann \cite{GangboMcCann96}, Uckelmann \cite{Uckelmann97},
McCann \cite{McCann99}, and Plakhov \cite{Plakhov04a}.

Imagine the periodic interval $X=Y = \R / 2\pi \Z=[0,2\pi[$ to parameterize a town built on
the boundary of a circular lake,  and let probability measures $\mu$ and $\nu$ represent
the distribution of students and available places in schools, respectively.
Suppose the distribution of students is smooth and non-vanishing but peaks sharply at the
northern end of the lake,  and the distribution of schools is smooth and non-vanishing
but peaks sharply at the southern end of the lake.
If the cost of transporting a student residing at location $\theta \in [0,2\pi]$
to school at location $\phi \in [0,2\pi]$ is presumed to be given in terms of the
angle commuted by $c(\theta,\phi) = 1 - \cos(\theta-\phi)$,  the most effective pairing of
students with places in schools is given by the measure in $\Gamma(\mu,\nu)$
which attains the minimum:
\begin{equation}\label{MKP}
\min_{\gamma \in \Gamma(\mu,\nu)} \int_{X \times Y} c(\theta,\phi) \; d\gamma(\theta,\phi).
\end{equation}

\begin{figure}[h]
\psfragscanon
\centering
\psfrag{o}{$o$}
\epsfig{file=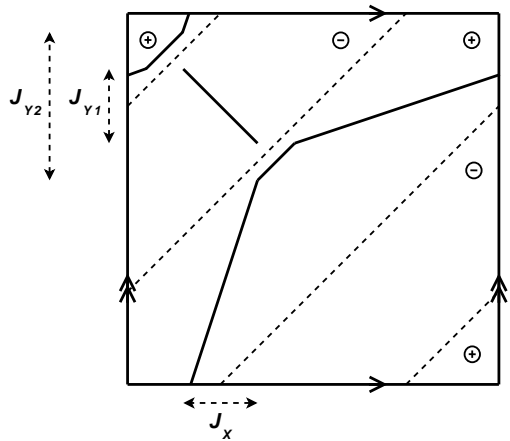, height=6cm}
\caption[An extreme doubly stochastic measure]
{\label{fig.extreme_doubly_stochastic}  \centering Schematic support of
the optimal measure from the example}
\end{figure}

\begin{figure}[h]
\psfragscanon
\centering
\psfrag{o}{$o$}
\epsfig{file=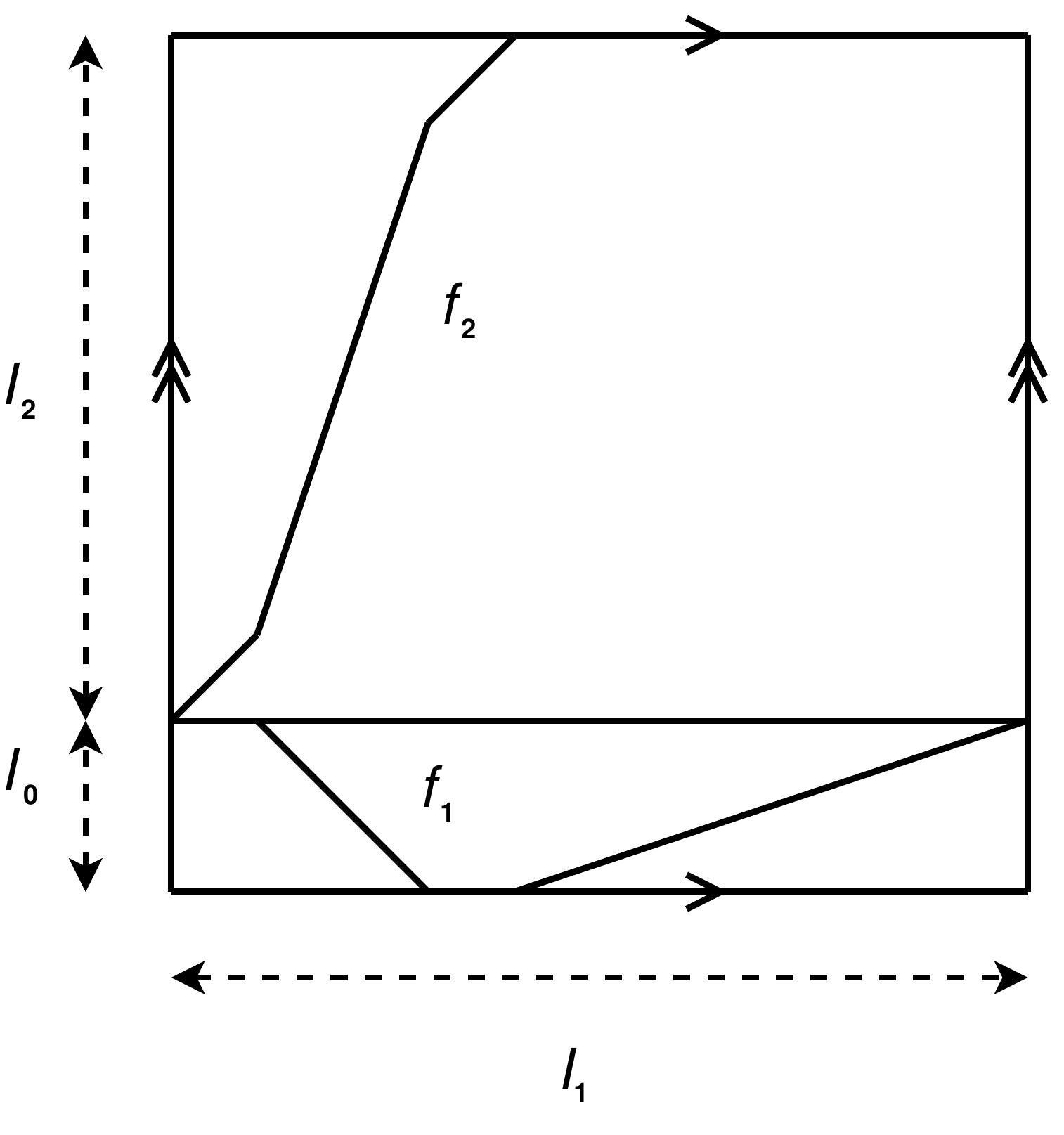, height=6cm}
\epsfig{file=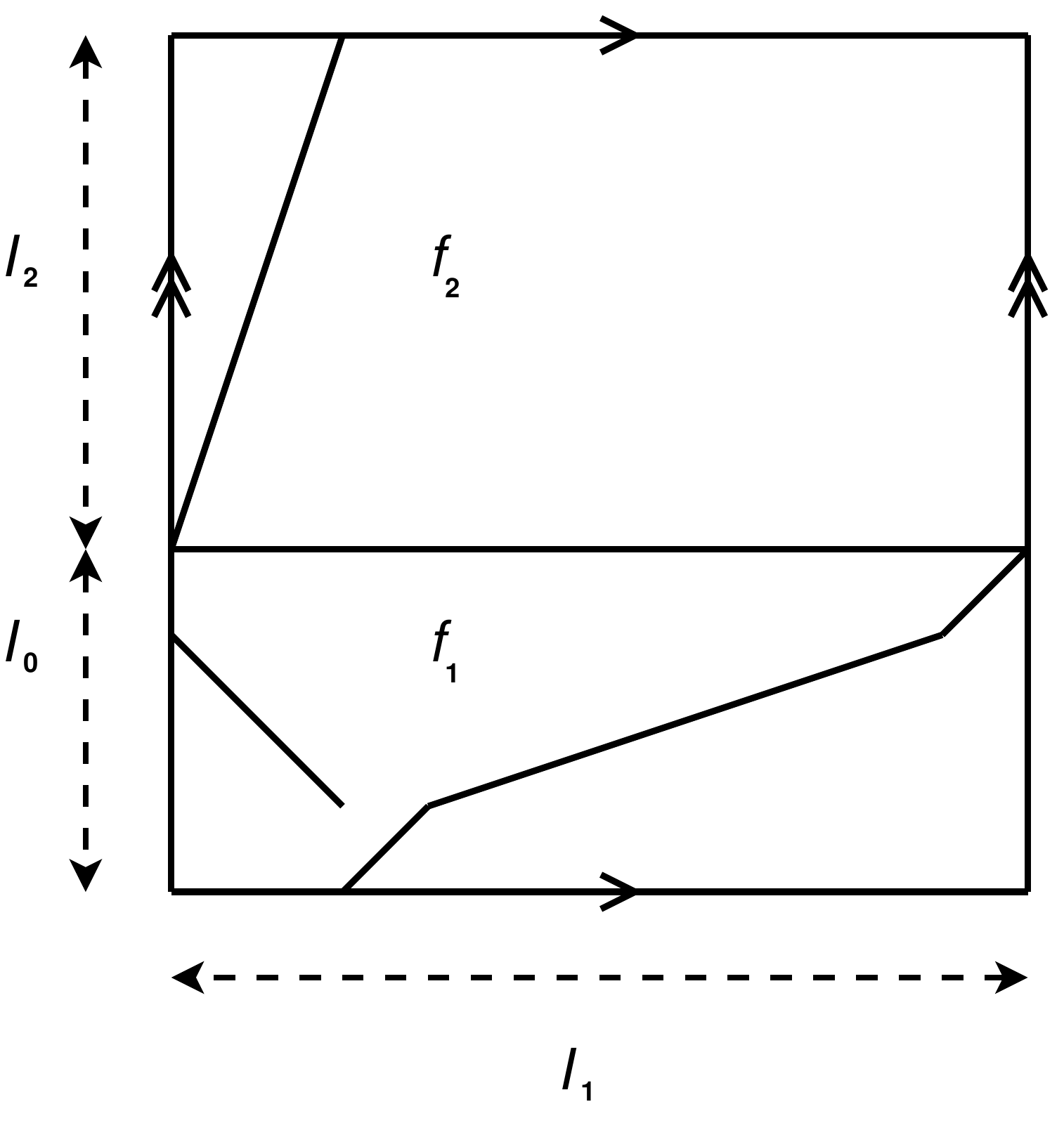, height=6cm}
\caption[Numbered limb system]
{\label{fig.example_numbered_limb}  \centering Two different numbered limb systems which represent Figure \ref{fig.extreme_doubly_stochastic}.}
\end{figure}

According to results of Gangbo and McCann \cite{GangboMcCann00},
this minimizer is unique, and its support is contained in the union of
the graphs of two maps ${\mathbf{t^\pm}}: X \longrightarrow Y$. A schematic
illustration is
given in Figure \ref{fig.extreme_doubly_stochastic}, where the restriction of
the support to the subsets marked by $\pm$ on the flat torus $X \times Y$ represent
$graph(\mathbf{t^+})$ and  $graph(\mathbf{t^-})$ respectively. The dotted lines mark
$\phi - \theta = \pm \frac{\pi}{2}, \pm \frac{3\pi}{2}$.
The necessary positivity of $\gamma[ J_{X} \times J_{Y1}]>0$
in this picture may be explained by observing that although
it is cost-effective for all students to attend a school where they live,  this is
incompatible with the concentration of students at the north end of the lake,  and
of schools at the south end.  Once this imbalance is corrected by sending a sufficient
number of northern students to southern schools by the map $\mathbf{t^-}$,
the remaining students can be assigned to school near their home using the map
$\mathbf{t^+}$.
Periodicity of graphs on the flat torus can be used to represent the support as a
numbered limb system in more than one way; see
Figure \ref{fig.example_numbered_limb}, which exploits the fact that the support
of $\gamma$ in Figure \ref{fig.extreme_doubly_stochastic} intersects
$X \times J_{Y2}$ in a graph and $X \times \left(Y - J_{Y1}\right)$
in an anti-graph.

Chiappori, Nesheim and McCann \cite{ChiapporiMcCannNesheim08p}
called the uniqueness hypothesis limiting the number of critical points
to at most one maximum and and at most one minimum in (\ref{subtwist})
the {\em subtwist} condition.
{Although it is satisfied in the example above,  it is an unfortunate fact
that the subtwist condition cannot be satisfied
by any smooth function $c(\theta,\phi)$ on a product of manifolds $X \times Y$
with more complicated Morse structures than the sphere.  It is an interesting
open problem to find a criterion on a smooth cost $c(\theta,\phi)$ on
$X=Y = \R^2/\Z^2$ which guarantees uniqueness of the minimum (\ref{MKP}) for all
smooth densities $\mu$ and $\nu$ on the torus.
Although we expect such costs to be generic,  not a single example of such a cost
is known to us.  Hestir and Williams criteria for extremality seems likely to remain
relevant to such questions, and it is natural to conjecture that the complexity of the
Morse structure of the manifold $X$ plays a role in determining the required number
of limbs in the system.}

\end{document}